# Every compact metric space that supports a positively expansive homeomorphism is finite

Ethan M. Coven[1] and Michael Keane[1]

*Wesleyan University, USA*

**Abstract:** We give a simple proof of the title.

A continuous map $f : X \to X$ of a compact metric space, with metric $d$, is called *positively expansive* if and only if there exists $\epsilon > 0$ such that if $x \neq y$, then $d(f^k(x), f^k(y)) > \epsilon$ for some $k \geq 0$. Here exponentiation denotes repeated composition, $f^2 = f \circ f$, etc. Any such $\epsilon > 0$ is called an *expansive constant*. The set of expansive constants depends on the metric, but the existence of an expansive constant does not.

Examples include all one-sided shifts and all expanding endomorphisms, e.g., the maps $z \mapsto z^n, n \neq 0, \pm 1$, of the unit circle. None of these maps is one-to-one, and with good reason. It has been known for more than fifty years that every compact metric space that supports a positively expansive homeomorphism is finite. This was first proved by S. Schwartzman [8] in his 1952 Yale dissertation. The proof appears in [4, Theorem 10.30]. (A mistake in the first edition was corrected in the second edition.) Over the years it has been reproved with increasingly simpler proofs. See [5],[6],[7].

The purpose of this paper is to give another, even simpler, proof of this result. The idea behind our argument is not new. After discovering the proof given in this paper, the authors learned from W. Geller [3] that the basic idea of the proof had been discovered in the late 1980s by M. Boyle, W. Geller, and J. Propp. Their proof appears in [6]. The result also follows from the theorem at the end of Section 3 of [1]. The authors of that paper may have been unaware of this consequence of their theorem, which they called "a curiosity based on the techniques of this work." That the result follows from the theorem in [1] was recognized by B. F. Bryant and P. Walters [2].

Nonetheless, the fact that this result has a simple proof remains less well-known than it should be, and so publishing it in this Festschrift is appropriate.

**Theorem.** *Every compact metric space that supports a continuous, one-to-one, positively expansive map is finite.*

**Remark.** Our statement does not assume that the map is onto. Such maps are sometimes called "homeomorphisms into." We have written the proof so that the fact that $f$ is one-to-one is used only once.

[1]Department of Mathematics, Wesleyan University, Middletown, CT 06459, USA, e-mail: ecoven@wesleyan.edu; mkeane@wesleyan.edu

*Proof.* Let $X$ be a compact metric space with metric $d$, let $f$ be a continuous, one-to-one, positively expansive map of $X$ into (not necessarily onto) itself, and let $\epsilon > 0$ be an expansive constant. Consider the following condition

(∗) there exists $n \geq 1$ such that if $d(f^i(x), f^i(y)) \leq \epsilon$ for $i = 1, 2, \ldots, n$, then $d(x, y) \leq \epsilon$, too.

Suppose that (∗) is not true. Then for every $n \geq 1$, there exist $x_n, y_n$ such that $d(f^i(x_n), f^i(y_n)) \leq \epsilon$ for $i = 1, 2, \ldots, n$, but $d(x_n, y_n) > \epsilon$. Choose convergent subsequences $x_{n_k} \to x$ and $y_{n_k} \to y$. Then $x \neq y$ and for every $i \geq 1$, $d(f^i(x_{n_k}), f^i(y_{n_k})) \to d(f^i(x), f^i(y))$. Now $d(f^i(x_{n_k}), f^i(y_{n_k})) \leq \epsilon$ for $i = 1, 2, \ldots, n_k$, so $d(f^i(x), f^i(y)) \leq \epsilon$ for all $i \geq 1$. But $f$ is one-to-one, so $f(x) \neq f(y)$. This contradicts positive expansiveness (with $f(x)$ and $f(y)$ in place of $x$ and $y$), so the condition holds.

Fix $k \geq 0$, and apply (∗) consecutively for $j = k, k-1, \ldots, 1, 0$, with $f^j(x)$ and $f^j(y)$ in place of $x$ and $y$. We get

(∗∗) for every $k \geq 0$, if $d(f^i(x), f^i(y)) \leq \epsilon$ for $i = k+1, k+2, \ldots, k+n$, then $d(f^i(x), f^i(y)) \leq \epsilon$ for $i = 0, 1, \ldots, k$, too.

Now cover $X$ by finitely many, say $N$, open sets of the form

$$U(z) := \{x \in X : d(f^i(x), f^i(z)) \leq \epsilon/2 \text{ for } i = 1, 2, \ldots, n\},$$

where $n$ is as in (∗).

If $X$ contains more then $N$ points, consider a finite subset containing $N+1$ points. For every $k \geq 0$, there exist $x_k \neq y_k$ in this subset such that $f^k(x_k)$ and $f^k(y_k)$ lie in the same set $U(z)$. Then by (∗∗), $d(f^i(x_k), f^i(y_k)) \leq \epsilon$ for $i = 0, 1, \ldots, k+n$. Since there are only finitely many pairs of these $N+1$ points, there exist $x \neq y$ such that $d(f^i(x), f^i(y)) \leq \epsilon$ for infinitely many $i \geq 0$, and hence for all $i \geq 0$. This contradicts positive expansiveness. □

## References

[1] ADLER, R. L., KONHEIM, A. G., AND MCANDREW, M. H. (1965). Topological entropy. *Trans. Amer. Math. Soc.* **114**, 309–319. MR175106
[2] BRYANT, B. F. AND WALTERS, P. (1969) Asymptotic properties of expansive homeomorphisms, *Math. Systems Theory* **3**, 60–66. MR243505
[3] GELLER, W. Personal communication.
[4] GOTTSCHALK, W. H. AND HEDLUND, G. A. (1955). *Topological Dynamics.* Amer. Math. Soc. Colloq Publ. **36**, Providence, RI. MR0074810
[5] KEYNES, H. B. AND ROBERTSON, J. B. (1969). Generators for topological entropy and expansiveness. *Math. Systems Theory* **3**, 51–59. MR247031
[6] KING, J. L. (1990). A map with topological minimal self-joinings in the sense of del Junco. *Ergodic Theory Dynam. Systems* **10**, 4, 745–761. MR1091424
[7] RICHESON, D. AND WISEMAN, J. (2004). Positively expansive homeomorphisms of compact metric spaces, *Int. J. Math. Math. Sci.* **54**, 2907–2910. MR2145368
[8] SCHWARTZMAN, S. (1952). On transformation groups. Dissertation, Yale University.